\documentclass[11pt]{amsart}
\usepackage{graphicx}
\usepackage{color}
\usepackage{latexsym}
\usepackage{amscd}
\usepackage{a4wide}
\usepackage{epsfig}
\usepackage{float}
\usepackage{amsmath, amssymb, amsthm, amsfonts, amsgen} 
\usepackage[outdir=./]{epstopdf}

\graphicspath{ {./} }

\usepackage{esint}

\usepackage{enumitem}
\usepackage{cite}

\usepackage{color}

\renewcommand{\[}{\begin{equation*}}
\renewcommand{\]}{\end{equation*}}
\newcommand{\grad}{\ensuremath{\nabla}}
\newcommand{\p}{\ensuremath{\partial}}
\newcommand{\ve}{\ensuremath{\varepsilon}}
\newcommand{\vp}{\ensuremath{\varphi}}
\newcommand{\al}{\ensuremath{\alpha}}
\newcommand{\be}{\ensuremath{\beta}}
\newcommand{\ga}{\ensuremath{\gamma}}
\newcommand{\ta}{\ensuremath{\tau}}
\newcommand{\ka}{\ensuremath{\kappa}}

\renewcommand{\th}{\ensuremath{\theta}}

\renewcommand{\ss}{\ensuremath{\subset}}
\newcommand{\sm}{\ensuremath{\setminus}}
\newcommand{\ra}{\ensuremath{\rightarrow}}
\newcommand{\mt}{\ensuremath{\mapsto}}

\newcommand{\mb}[1]{\ensuremath{\;\mbox{#1}\;}}
\newcommand{\cd}{\ensuremath{\cdot}}
\newcommand{\ot}{\ensuremath{\otimes}}
\newcommand{\ti}{\ensuremath{\times}}

\newcommand{\Id}{\ensuremath{\mbox{Id}}}

\newcommand{\wt}[1]{\widetilde{#1}}
\newcommand{\ol}[1]{\overline{#1}}

\newcommand{\cof}{\textnormal{cof}\;}

\newcommand{\dor}{\ensuremath{\dot{r}}}
\newcommand{\doz}{\ensuremath{\dot{z}}}
\newcommand{\ddoz}{\ensuremath{\ddot{z}}}
\newcommand{\dod}{\ensuremath{\dot{d}}}

\newcommand{\ddor}{\ensuremath{\ddot{r}}}

\newcommand{\R}{\mathbb{R}}

\newcommand{\N}{\mathbb{N}}

\newcommand{\D}{\Delta}
\renewcommand{\d}{\delta}
\newcommand{\A}{\mathcal{A}}

\newcommand{\La}{\mathcal{L}}

\newtheorem{de}{Definition}[section]
\newtheorem{thm}[de]{Theorem}

\newtheorem{lem}[de]{Lemma}

\newtheorem{re}[de]{Remark}

\allowdisplaybreaks
\begin{document}
\title[On $C^1-$stationary points of a polyconvex functional and finite BOP-theory]{On the assembly of $C^1-$stationary points of a polyconvex functional and finite BOP-theory} 
\author[M Dengler]{M. Dengler}
\address[M. Dengler]{Fliederweg 1, 72189 Vöhringen, Germany.}
\email{m.dengler@web.de}

\keywords{Calculus of Variations, elasticity, polyconvexity, regularity, BOP-theory, ODE} 
\subjclass[2020]{49N60, 73C50, 34A34, 34B15}

\maketitle 
\begin{abstract}
\noindent In this work the following energy is considered
\begin{equation*}I(u)=\int\limits_B{\frac{1}{2}|\grad u|^2+\rho(\det\grad u)\;dx},\label{eq:SA.1.1}\end{equation*}
where $B\ss\R^2$ denotes the unit ball, $u\in W^{1,2}(B,\R^2)$ and $\rho:\R\ra\R_0^+$ smooth and convex with $\rho(s)=0$ for all $s\le0$ and $\rho$ becomes affine when $s$ exceeds some value $s_0>0.$ Additionally, we may impose $(M\in \N\sm\{0\})-$covering maps as boundary conditions in a suitable fashion.\\
For such situations we then construct radially symmetric $M-$covering stationary points of the energy, which are at least $C^1$ (in some circumstances even $C^\infty),$ and verify more refined properties, which these stationary points need to satisfy. We do so by following the strategy first and foremost developed by P. Bauman, N. C. Owen, and D. Phillips (BOP) confirming and generalising that the method remains valid beyond the $M=2-$case for an arbitrary $M.$ Furthermore, as far as we know, this is the first treatise of BOP-theory in finite elasticity. The finiteness, not imposing such strict conditions, allows for a richer class of possible behaviours of the stationary points, making it more difficult to completely determine them.
\end{abstract}

\section{Introduction}
Let $B\ss\R^2$ be the unit ball and define the functional $I:W^{1,2}(B,\R^2)\rightarrow \R$ by
\begin{equation}
I(u):=\int\limits_{B}{\frac{1}{2}|\grad u|^2+\rho(\det \grad u)\; dx}
\label{eq:1.1}
\end{equation}
for all $u\in W^{1,2}(B,\R^2).$ The function $\rho\in C^{\infty}(\R)$ is defined by
\begin{equation}
\rho(s)=\left\{\begin{array}{ccc}
0& {\mbox{if}}& s\le0,\\
\rho_1(s)& {\mbox{if}}& 0\le s\le s_0,\\
\gamma s+\ka&{\mbox{if}}& s_0\le s,
\end{array}
\right.
\label{eq:1.2}
\end{equation}
for some constants $\gamma>0,$ $s_0\ge0$ and $\ka\ge-\ga s_0.$ 
Here $\rho_1:[0,s_0]\rightarrow \R$ is a smooth and convex function on $[0,s_0]$ satisfying the boundary conditions $\rho_1(0)=0$ and $\rho_1(s_0)=\gamma s_0+\ka$ and the connections need to be in such a way that $\rho$ is smooth everywhere. Note that $\rho$ is convex on the whole real line. Hence, the complete integrand is polyconvex. Recall, again, we call $W:\R^{2\ti2}\ra\ol{\R}$ polyconvex, if there exists a convex function $g:\R^5\ra\ol{\R}$ s.t. $W(\xi)=g(\xi,\det\xi)$ for all $\xi\in\R^{2\ti2}.$
The behaviour of the functional depends mainly on the parameter $\ga.$ If $\ga\ra0$ then the functional turns into the well known Dirichlet energy. In the regime $0<\ga<1$ the functional is uniformly convex. If $\ga\ge1$ then the functional is genuinely polyconvex.\\

Now we define the $M-$covering map $u_M:B\rightarrow B$ for $M\in\N,$ $M\ge1$ via its representative
\begin{eqnarray}
u_M:[0,1]\ti[0,2\pi)&\rightarrow& B,\\
(R,\th)&\mapsto& R(\cos M\th,\sin M\th).
\end{eqnarray}
Furthermore, we will use the notation $e_{MR}:=e_R(M\th):=(\cos M\th,\sin M\th)^T$ and $e_{M\th}:=e_\th(M\th):=(-\sin M\th,\cos M\th)^T$ for all $M\in\N\sm\{0\}$ and we introduce the set of admissible functions
\[\A_{u_M}:=\{u\in W^{1,2}(B,\R^2): u=u_M \;\mb{on}\; \p B\},\]
where the boundary condition must be understood in the trace sense. Moreover, $\A_{u_M}\not=\emptyset$ and the functional $I$ attains its minimum in $\A_{u_M}.$ \\

\textbf{Notation:} For a $2\times2-$matrix $A$ we denote the determinant by $d_A:=\det A$ (where we drop the $A$ for simplicity)  and the cofactor is given by 
\begin{equation}
\cof A=\begin{pmatrix}a_{22} &-a_{21}\\-a_{12} & a_{11}\end{pmatrix}.
\label{eq:1.3}
\end{equation}
For two vectors $a\in \R^n,b\in \R^m$ we define the tensor product $a\ot b\in \R^{n\ti m}$ by $(a\ot b)_{i,j}:=(ab^T)_{i,j}=a_ib_j$ for all $1\le i\le n,$ $1\le j\le m.$\\

As usual we are interested in the behaviour of the minimizers or more generally the class of stationary points:
\begin{de}[stationary point]
We call $u\in\A_{u_M}$ a stationary point of the functional \eqref{eq:1.1} if it satisfies the Euler-Lagrange equation (ELE) in the weak form, which is given by
\begin{equation}
\int\limits_{B}{(\grad u+\rho'(d_{\grad u})\cof \grad u) \cdot\grad\vp\; dx}=0\; {\mbox{for all}}\; \vp\in W_0^{1,2}(B,\R^2).
\label{eq:ELE1.1}
\end{equation}
\end{de}
\vspace{0.3cm}

The set of radially symmetric $M-$covering maps
\[\A_r^M=\{u\in W^{1,2}(B, \R^2)|\exists r:[0,1]\rightarrow\R\; \mbox{s.t.}\; u(x)=r(R)e_R(M\theta)\;\mbox{and}\;r(1)=1\}\ss \A_{u_M}\]
 will play a key-role in this paper. Again, $\A_r^M\not=\emptyset$ and the functional $I$ attains its minimum in $\A_r^M,$ as well. Now restricting the class of test functions to $M-$covering radial symmetric maps, too, then the ELE becomes a boundary value problem (BVP) for $r:$ 

\begin{de}[BVP] For any $M\in \N\sm\{0\}$ and $\ga\in (0,\infty)$ we call $r\in L^2((0,1),R^{-1}dR)\wedge \dor \in L^2((0,1),RdR)$ a distributional solution of the BVP if it satisfies
\begin{align}\left\{\begin{array}{ccc}
\frac{M^2r}{R}-\dor-R\ddor=M\rho''(d)\dod r & {\mbox{in}}& \mathcal{D}'(0,1),\\
r(0)=0,&r(1)=1.&
\end{array}\label{eq:0.0.1}
\right.\end{align}\\
\end{de}

Furthermore, we will be interested in the concrete shapes of $\rho$ and $r$ and therefore we will need the following notions:

\begin{de}[lift-off functions]
Let $h\in C(\R)$ be a real valued function.\\

\begin{enumerate}
\item We call $h$ an immediate lift-off function, if $h(0)=0$ and $h(R)>0$ for any $R>0.$ 
\item We call $h$ a delayed lift-off function, if there is $\d>0$ s.t. $h\equiv0$ on $[0,\d]$ and $h(R)>0$ for any $R>\delta.$

\end{enumerate}
\end{de}\vspace{0.3cm}

This brings us in the position to give the central statement:
 
\begin{thm} Let $M\in \N\sm\{0\},$ $\ga\in (0,\infty),$ let $I$ be the energy as given in \eqref{eq:1.1}, and let $\rho$ as defined in \eqref{eq:1.2} be arbitrary. \\

Then the following statements are true:
\begin{enumerate}
\item Every map $u=re_{MR}\in \A_r^M$ with $r$ satisfying the BVP \eqref{eq:0.0.1} distributionally is a stationary point of $I$ in the full class $\A_{u_M}.$
\item If $M=1$ then $u=\Id$ is the unique global minimizer of $I$ in the full class $\A_{u_M}.$
\item If $M\ge2$ then every map $u=re_{MR}$ with $r$ satisfying the BVP \eqref{eq:0.0.1} distributionally possesses the following properties:
\begin{enumerate}
\item $u\in C^1(\ol{B},\R^2).$
\item $r\in C^1([0,1])\cap C^\infty((0,1])$ and $r,\dor\ge0$ for all $R\in[0,1]$ and $r(0)=\dor(0)=0.$ 
\item $d:=\det \grad u=\frac{Mr\dor}{R}\in C([0,1])\cap C^\infty((0,1])$ and $d\ge0$ in $[0,1]$ with $d(0)=0$ and $\dod\ge0$ in $(0,1].$
\item If $\rho$ is lifting-off delayed then $r,d\in C^\infty([0,1])$ and $u=re_{MR}\in C^\infty(\ol{B},\R^2).$
\item If $\rho$ is lifting-off immediately then $r$ is either a delayed lift-off solution with $r\in C^\infty([0,1])$ or 
an immediate lift-off function with $r\in C^1([0,1])\cap C^\infty((0,1]),$ but not necessarily any better, and $r(0)=\dor(0)=0.$ 
If it is additionally assumed that $u=re_{MR}\in W^{2,2}(\ol{B},\R^2)$ then $r,d\in C^\infty([0,1])$ and $u=re_{MR}\in C^\infty(\ol{B},\R^2).$
\end{enumerate}
\end{enumerate}
\label{thm:1.3 MT}
\end{thm}\vspace{0.5cm}

\begin{re}
Notice, that in theorem \ref{thm:1.3 MT}.(1) it is shown that the involved maps are stationary points of the functional $I$ in the full class $\A_{u_M}$ however, for $M\ge2$ it remains open, if theses constructed maps actually are the global minimizers of the energy $I$.
\end{re}\vspace{3mm}

What does the literature tell us about the regularity of stationary points/minimizers of the energy $I$? This functional has recently been introduced in \cite{D2}. There it is shown that for arbitrary $\ga\in (0,\infty)$ any stationary point of the energy $I$ subject to arbitrary $L^2-$ boundary data needs to be locally Hölder-continuous and a higher-order regularity result is obtained. Moreover, it is shown that if $I$ is uniformly elliptic (that is $\ga\in (0,1)$) any stationary point of the energy $I$ must be locally smooth. 
Furthermore, the classical result by Acerbi and Fusco \cite{AF87} guarantees that any minimizer/stationary point of $I$ under suitable boundary conditions has to be smooth up to a nullset. Although, the fairly symmetric stationary points constructed here are at least of class $C^1,$ less symmetric stationary points with low regularity could still occur: Müller and Sverak \cite{MS03}, Kristensen and Taheri \cite{KT03}, and Szekelyhidi, Jr. \cite{Sz04}  constructed, via Gromov's convex integration method, smooth poly- and quasiconvex integrands s.t. the stationary point/ local minimizer is everywhere Lipschitz but nowhere $C^1$.
Elastic situations subject to double (or $M-$) covering boundary data have been intensely discussed in \cite{BallOP,JB14,BeDe21,D1,BOP91MS}. Additional everywhere regularity results can be found in \cite{FH94, C14, CLM17,BOP92, K11,K13} and for a list of partial regularity results regarding poly- or quasiconvex integrands see \cite{FH94,EM01,CY07,Ev86,AF87,KT03, SS09,CG20,KM07,EVGA99}.\\

In this paper we follow the method devised initially by P. Bauman, N.C. Owen and D. Phillips (BOP) in two striking papers \cite{BOP91, BOP91MS} and extended by Yan \cite{Y06} and Yan and Bevan \cite{BY07}. Indeed, BOP consider a infinite nonlinear elastic situation, where the integrand $W$ depends on $d$ in such a way, that $W(d)\ra+\infty$ if $d\ra0^+$ or $d\ra+\infty$ and $W(d)=+\infty$ if $d\le 0$ and they obtain a higher-order regularity result showing that it is enough for an equilibrium solution to be of class $C^{1,\be}$ for some $\be\in(0,1]$ for it to be already fully smooth. In \cite{Y06} it is then shown that the previous assumption can be relaxed to $C^{1}\cap W^{2,2}$ and that this is optimal in a sense. Additionally, BOP construct a radial symmetric double covering singular equilibrium solution, which is $C^1\sm C^{1,\be}$ for any $\be\in(0,1].$ Bevan and Yan \cite{BY07} then show, that this map is the unique minimizer in a substantial part of the admissible set. The latter assembly of the counterexample is important since, we follow this part very closely. For more regularity results in infinite nonlinear elasticity  see \cite{JB17,BK19,FR95}. \\

\textbf{Plan for the paper:} In §.2 we start by developing the classical BOP-theory. The general goal of §.2 is to get an understanding of the regularity of $r$ and $d.$ Initially, we discuss basic results of $r$ in lemma \ref{lem:2.1}. In lemma \ref{lem:Mdet} we obtain regularity and non-negativity for $d$ and $\dod,$ and then analogous results for $r$ and $\dor$ in lemma \ref{lem:2.3}. This is followed by two technical lemmas describing the behaviour of two important auxiliary functions $z$ and $f.$ These are then used to obtain lemma \ref{lem:3.2.100}, which is the central statement of §.2 discussing the notions of delayed and immediate lift-off solutions and their corresponding regularity. Lastly, lemma \ref{lem:2.7} shows point $(1)$ of theorem \ref{thm:1.3 MT}. In §.3 we then explore the solutions, starting-off, following the newer theory developed by Yan and Bevan, by discussing, in lemma \ref{lem:3.1}, limit taking at the origin involving $\ddor.$ In §.3.1 for an arbitrary $\rho$ delayed lift-off solutions are investigated. Finally, §.3.2 and §.3.3 then distinguish between situations, where $\rho$ is lifting-off delayed, and such, where $\rho$ is  immediately lifting-off and the corresponding solutions are analysed once again.

\section{Classical BOP-theory}
In this paragraph it is shown that functions are at least of class $C^1$ on the whole interval $[0,1]$. Again, we use the method invented in \cite{BOP91, BOP91MS}. In particular, this section will follow very closely the latter one without us mentioning it all the time.\\

We start our discussion by recalling the ELE as given in \eqref{eq:ELE1.1}. Plugging $u=re_{MR}\in \A_r^M$ and test functions of the form $\vp(x)=ge_{NR} \;\mbox{where}\; g\in C_c^\infty((0,1)), N\in \N\sm\{0\}$ into \eqref{eq:ELE1.1} yields the BVP \eqref{eq:0.0.1}.
The first statement makes this rigorous and establishes first simple results of the radial part $r.$

\begin{lem}[Elementary properties]Let $M\in\N,$ $M\ge1$ and $u\in \A_r^M.$
\begin{enumerate}[label=(\roman*)]
\item Then $u\in \A_r^M$ if and only if $r$ is absolutely continuous on each compact subset of $(0,1]$ and $r\in L^2((0,1), R^{-1}dR)$ and $\dor\in L^2((0,1), RdR).$\\ 
\item If $M\not=N$ then $W$ is a Null-Lagrangian in the class of $\A_r^M$ and $N-$covering test functions, i.e.
\[\int\limits_B{(\grad u+\rho'(d)\cof\grad u)\cd\grad\vp\;dx=0 }\]
holds for all $u\in\A_r^M$ and all test functions of the form $\vp(x)=g(R)e_R(N\th)$ where $g\in C_c^\infty((0,1)),$ if $M\not=N.$
\item Assume, additionally, $u\in \A_r^M$ solves
\[\int\limits_B{(\grad u+\rho'(d)\cof\grad u)\cd\grad\vp\;dx=0 }\] weakly, for any test functions of the form $\vp(x)=g(R)e_R(N\th)$ where $g\in C_c^\infty((0,1)).$ Then $r$ satisfies the ODE
\begin{equation}
\left(\frac{M^2r}{R}+M\rho'(d)\dot{r}\right)=\left(R\dot{r}+M\rho'(d)r\right)^{\cdot} \;\mbox{in}\; D'((0,1)).
\label{eq:3.15}
\end{equation}
\item Assume $r$ solves \eqref{eq:3.15}. Then $r\in C([0,1])$ and $r(0)=0.$
\item Moreover, $r\in C^{\infty}((0,1]).$
\end{enumerate}
\label{lem:2.1}
\end{lem}
\begin{proof}
(i): It is straightforward to see that for every $u\in\A_r^M$ it holds that
\[\|u\|_{W^{1,2}}^2=\int\limits_B{|u|^2+|\grad u|^2\;dx}=2\pi\int\limits_0^1{r^2R+\frac{M^2r^2}{R}+\dor^2R\;dR}\]
Since the LHS is finite the RHS needs to be finite, too, implying $r\in L^2((0,1), R^{-1}dR)$ and $\dor\in L^2((0,1), RdR).$\\
Let $[a,b]\ss(0,1].$
Since $r\in L^2((0,1], R^{-1}dR)$ it follows
\begin{equation}
\int\limits_a^b{r^2\;dR}\le b\int\limits_a^b{r^2\;\frac{dR}{R}}<\infty. 
\label{eq:3.24}
\end{equation}
Hence, $r\in L_{loc}^2((0,1]).$ Similar $\dor\in L^2((0,1], RdR)$ leads to 
\begin{equation}
\int\limits_a^b{\dor^2\;dR}\le\frac{1}{a}\int\limits_a^b{\dor^2\;RdR}<\infty
\label{eq:3.25}
\end{equation}
which implies $\dor\in L_{loc}^2((0,1])$ and $r\in W_{loc}^{1,2}((0,1]).$ Then $r$ agrees up to a set of measure zero with a function $\tilde{r}$ on $(0,1],$ where $\tilde{r}$ is absolutely continuous on any compact subset of $(0,1]$ (as always we identify $r$ with $\tilde{r}$). The latter is a consequence of \eqref{eq:3.25} and the fundamental theorem of calculus for Sobolev functions, see \cite[U1.6, p.71-72]{A12}.\\

(ii)-(iii): First we calculate some important quantities:
\begin{align*}
\grad u=&{u,}_{R}\ot e_R+u,_{\ta}\ot e_\th=\dot{r}e_{MR}\ot e_R+\frac{Mr}{R}e_{M\th}\ot e_\th&\\
\grad \vp=&\dot{g}e_{NR}\ot e_R+\frac{Ng}{R}e_{N\th}\ot e_\th&\\
\cof\grad u=&\frac{Mr}{R}e_{MR}\ot e_R+\dot{r}e_{M\th}\ot e_\th&\\
d=&\det\grad u=\frac{1}{2}\grad u\cdot \cof\grad u=\frac{Mr\dot{r}}{R}&\\
\grad u\cdot\grad \vp=&\dot{g}\dot{r}e_{MR}\cdot e_{NR}+\frac{MNgr}{R^2}e_{M\th}\cdot e_{N\th}&\\
\cof\grad u\cdot\grad \vp=&\frac{M\dot{g}r}{R}e_{MR}\cdot e_{NR}+\frac{Ng\dot{r}}{R}e_{M\th}\cdot e_{N\th}&\\
e_{MR}\cdot e_{NR}=&e_{M\th}\cdot e_{N\th}=\cos(M\th)\cos(N\th)+\sin(M\th)\sin(N\th)=\cos((M-N)\th).&
\end{align*}

 Hence, the ELE \eqref{eq:ELE1.1} becomes
\begin{equation}
\int\limits_{0}^1\int\limits_{0}^{2\pi}{\left(\dot{g}\dot{r}+\frac{Mr\dot{g}}{R}\rho'(d)+\frac{MNgr}{R^2}+\rho'(d)\frac{Ng\dot{r}}{R}\right)\cos((M-N)\th)\;d\th RdR}=0
\label{eq:3.21}
\end{equation}
which is automatically true for all $r$ if $M\not=N$ and for $M=N$ takes the form
\begin{equation}
2\pi\int\limits_{0}^1{\dot{g}\left(\dot{r}R+M\rho'(d)r\right)+g\left(\frac{M^2r}{R}+M\rho'(d)\dot{r}\right)\;dR}=0
\label{eq:3.22}
\end{equation}
for all $g\in C_c^{\infty}((0,1)).$\\

(iv): Now we show that for any $r$ with $r\in L^2((0,1), R^{-1}dR)$ and $\dor\in L^2((0,1), RdR),$ it must hold that $\lim\limits_{R\ra0}r(R)=0.$ Suppose not, then wlog.\!\! there exists a strictly montonic decreasing sequence $\{R_j\}_{j\in \N}$ s.t. $R_j\ra0$ for $j\ra\infty$ and $|r(R_j)|>2\ve$ for any $j\in\N.$
Then since $\dor\in L^2((0,1), RdR),$ we can find $N\in\N$ so large that for any $n\ge N$ it holds
$\int\limits_0^{R_n}\dor^2R\;dR<\ve^2.$ Using the latter together with the fundamental theorem of calculus, and H\"older's inequality, then for any $n\ge N$ and any $R\in[R_n/e_n,R_n]$ with $e_n:=1-1/n$ (Note, that wlog. we can assume $R_{n+1}\le \frac{R_n}{e_n}$ if not consider the sequence $\wt{R}_{n+1}:=\min\{R_{n+1},\frac{R_n}{e_n}\}$) we obtain
 \[|r(R)-r(R_n)|=\left|\int_{R}^{R_n}\dot r(R')\, dR'\right|\le 
 \left| \int_{R}^{R_n}\frac1{R'}\, dR'\right|^\frac12\cdot 
 \left| \int_{R}^{R_n} 
 \left|\dot r(R')\right|^2R'\, dR'
 \right|^\frac12\le |\ln e_n|^\frac12 \ve.
 \]      
 Therefore by the reverse triangle inequality we have $|r(R)|>\ve$ for any $R\in[R_n/e_n,R_n]$. Then
 \[
 \int_{0}^{1}{r^2\over R}\,d R\ge\sum\limits_{n=N}^\infty\int_{R_n/e_n}^{R_n}{r^2\over R}\,d R> \sum\limits_{n=N}^\infty\int_{R_n/e_n}^{R_n}{\ve^2\over R}\,d R=\ve^2\sum\limits_{n=N}^\infty\left|\ln \left(1-\frac1{n}\right)\right|=+\infty,
 \]
 contradicting the integrability of ${r^2\over R}$ and showing that $r\in C([0,1])$ with $r(0):=0.$\\

(v): Next we show that $\dor\in C((0,1]).$\\
Let
\begin{equation}
q(w,a,b)=bw+M\rho'\left(\frac{Mwa}{b}\right)a
\label{eq:3.23}
\end{equation}
with $1\ge b>0,a,w\in\R.$
Then $q(\cdot,a,b)$ is a homeomorphism from $\R$ to $\R$ for all $a\in \R$ and  $b\in(0,1].$ Indeed, for $a=0,$ $q(w,0,b)=bw$ is a homeomorphism from $\R$ to $\R.$ Now let $b\in(0,1]$ and $a>0.$ Then $w\mt bw$ is strictly monotonically increasing and maps $\R$ to $\R$ continuously. Moreover, $w\mt M\rho'\left(\frac{Mwa}{b}\right)a$ is monotonically increasing and continuous as well, hence, $w\mt q(w,a,b)$ is a homeomorphism on $\R.$ If $a<0$ then $w\mt\rho'\left(\frac{Mwa}{b}\right)$ decreases but $w\mt M\rho'\left(\frac{Mwa}{b}\right)a$ still increases, and we can argue as above.\\
Assume that there exists $R_0\in(0,1]$ and a sequence $R_j\ra R_0$ for $j\ra\infty,$ s.t. $\dor(R_j)\ra\pm\infty.$ But then
\[q(\dor(R_j),r(R_j),R_j)=R_j\dor(R_j)+M\rho'\left(\frac{M\dor(R_j)r(R_j)}{R_j}\right)r(R_j)\ra\pm\infty+c=\pm\infty\]
which is impossible, since $R\ra q(\dor(R),r(R),R)$ is continuous on $(0,1].$ Hence, $\dor\in C((0,1]).$\\
As a last step we improve the regularity to $C^\infty((0,1]).$\\
Using the ODE we can represent $q$ by
\begin{equation}
q(\dor(R),r(R),R)=R\dor(R)+M\rho'(d)r(R)=c-\int\limits_R^1{\frac{M^2r(R')}{R}+M\rho'(d)\dor(R')\;dR'}
\label{eq:3.23}
\end{equation}
Then by the analysis above, $q$ is actually $C^1((0,1]).$ Further the derivative w.r.t. $w$ is $\p_wq=b+M^2\rho''(\frac{Mwa}{b})\frac{a^2}{b}>0$ for all $1\ge b>0,a,w\in\R.$ Then the implicit function theorem gives full regularity $r\in C^\infty((0,1]).$\end{proof}\vspace{0.5cm}

Collecting these results, from now on, we will consider solutions $r\in C([0,1])\cap C^\infty((0,1])$ to the boundary value problem
\begin{align}\left\{\begin{array}{ccc}
Lr=M\rho''(d)\dod r & {\mbox{in}}& (0,1),\\
r(0)=0,&r(1)=1,&
\end{array}\label{eq:3.0.1}
\right.\end{align}
where $L$ refers to the linear part of the considered ODE, i.e.
\[Lr(R):=\frac{M^2r}{R}-\dor-R\ddor\]
for all $R\in(0,1).$ \\

We keep following \cite{BOP91MS} so it might be good to outline the strategy: Firstly we will consider three auxiliary functions, namely  the determinant $d$, $z,$ and $f$ as defined below, which depend on $r$ and $\dor.$ Studying their behaviour, in particular, close to the origin will then reduce the number of possibilities how $d,$ $z,$ and $f$ can behave. Then we can discuss these cases one-by-one and finally determine, among other things, the regularity of $r$ and $\dor$. \\

We start with the following observations on the determinant. 
 
\label{sec:3.2}
\begin{lem}Let $M\in\N,\; M\ge1$ and assume that $r$ solves the BVP \eqref{eq:3.0.1}. Then $d\in C^\infty((0,1])$ and $\dod\ge0$ in $(0,1].$ Moreover, it holds that $d\in C([0,1])$ and $d\ge0$ in $[0,1].$
\label{lem:Mdet}
\end{lem}
\begin{proof}
The smoothness of $d$ in $(0,1]$ follows by the smoothness of $r$ in the same interval.
\[
\dod=\left(\frac{Mr\dor}{R}\right)^\cdot=\frac{M\dor^2}{R}+\frac{Mr\ddor}{R}-\frac{Mr\dor}{R^2}
\label{eq:3.d1}
\]
Now by multiplying the strong form of the ODE \eqref{eq:3.15} by $\frac{r}{R}$ we can express the term $r\ddor$ as
\begin{eqnarray}
r\ddor=\frac{M^2r^2}{R^2}-\frac{r\dor}{R}-M\rho''(d)\frac{r^2}{R}\dod.
\label{eq:3.d2}
\end{eqnarray}
Hence,
\[\dod=\frac{M\dor^2}{R}-\frac{Mr\dor}{R^2}+\frac{M^3r^2}{R^3}-\frac{Mr\dor}{R^2}-M^2\rho''(d)\frac{r^2}{R^2}\dod\]
and rearranging the equation yields, 
\begin{eqnarray}
\dod&=&(1+M^2\rho''(d)\frac{r^2}{R^2})^{-1}\left[\frac{M\dor^2}{R}-\frac{2Mr\dor}{R^2}+\frac{M^3r^2}{R^3}\right] \nonumber\\
&=&\frac{M[(R\dor-r)^2+(M^2-1)r^2]}{R^3+M^2\rho''(d)r^2R}\ge0.
\label{eq:3.d3}
\end{eqnarray}

Assume now that $\lim\limits_{R\ra0}d(R)\in[-\infty,0).$ By the smoothness of $d$ in $(0,1],$ there exists $\d>0$ s.t. $d(R)<0$ for all $R\in(0,\d).$ This implies that for all $R\in(0,\d)$ either $r(R)>0$ and $\dor(R)<0$ or $r(R)<0$ and $\dor(R)>0.$ Consider the first case $r,-\dor>0$ on $(0,\d).$ By the mean value theorem we get that there exists a $\xi\in(0,\frac{\d}{2})$ s.t. $\dor(\xi)=\frac{2r(\frac{\d}{2})}{\d}>0,$ contradicting $\dor(R)<0$ for all $R\in(0,\d).$ (Analogously, for the other case). This proves $\lim\limits_{R\ra0}d(R)\in[0,+\infty].$ By $\dod\ge0$ and $r$ smooth in $(0,1]$ $d$ cannot attain $+\infty.$ Therefore, the limit $d(0):=\lim\limits_{R\ra0}d(R)\in[0,\infty)$ exists and is nonnegative. Again by $\dod\ge0,$ $d$ remains nonnegative throughout the whole interval $[0,1].$  \end{proof}\vspace{0.5cm}

As a consequence, the non-negativity and the monotonic growth of $d$ are transferred on to $r.$

\begin{lem} Let $M\in\N,$ $M\ge1$ and $r$ solves the BVP \eqref{eq:3.0.1}. Then $r(R)\ge0$ for all $R\in[0,1]$ and $\dor(R)\ge0$ for all $R\in(0,1].$ 
\label{lem:2.3}
\end{lem}
\begin{proof}
By Lemma \ref{lem:Mdet} we know that
\[d=\frac{Mr\dor}{R}=\frac{M(r^2)^\cdot}{2R}\ge0.\]
Hence, $(r^2)^\cdot\ge0$ for all $R\in(0,1].$ For the sake of a contradiction, assume that there exists $R_0\in(0,1]$ s.t. $r(R_0)<0.$ Then by the continuity of $r$ and since $r^2$ grows monotonically, $r$ remains negative up to the boundary, i.e. $r(R)\le r(R_0)<0$ for all $R\in[R_0,1].$ This is not compatible with the boundary condition $r(1)=1$ yielding $r\ge0$ in $[0,1].$\\
The second claim follows in a similar fashion. Again we make the assumption that there exists $R_0\in(0,1]$ s.t. $\dor(R_0)<0.$ By continuity, there even exists an interval $(R_1,R_2]$ with $0<R_1<R_2\le1$ s.t. $\dor(R)<0$ for all $R\in (R_1,R_2]$ then by monotonicity of $r^2$ we know that $0\le r(R_2)-r(R_1)$ but on the other hand by the fundamental theorem of calculus we have
\[0\le r(R_2)-r(R_1)=\int\limits_{R_1}^{R_2}{\dor(R)\;dR}<0\]
leading again to a contradiction. Hence, $\dor\ge0$ in $(0,1].$\end{proof}\vspace{0.5cm}

Next we introduce the function $z(x):=\frac{1}{2}|\grad u(x)|^2+f(\det\grad u(x))$ for all $x\in\ol{B},$ where $f(d):=d\rho'(d)-\rho(d)$ for all $d\in\R.$ In the following lemma it is shown that $z$ satisfies a maximum principle in $\ol{B}\sm\{0\}$. This follows closely \cite[thm $3.2-3.3$]{BOP91}.

\begin{lem} Let $M\in\N,\; M\ge2$ and assume $r$ solves the BVP \eqref{eq:3.0.1}. Then $z$ satisfies the strong maximum principle in $(0,1].$
\label{lem:3.4}
\end{lem}
\begin{proof}
It is enough to show that $z$ is a subsolution to an elliptic equation, i.e. $\D z+c_M\rho''(d)\doz\ge0$ in $(0,1],$ where $c_M:=\frac{M}{2(M-1)}.$ This is indeed enough to apply the strong maximum principle, see \cite[\S6.4.2 Thm 3]{LE10}. \\
Initially, note that for $u\in\A_r^M,$  $z(x)=z(R),$ $z(R)=\frac{\dor^2}{2}+\frac{M^2r^2}{2R^2}+f(d),$ $f(d)=d\rho'(d)-\rho(d)$ and $\D z=\frac{\doz}{R}+\ddoz.$\\
Now in order to calculate $\D z$ we first need to calculate $\doz$ and
 $\ddoz.$ Taking the derivative of $z$ wrt. $R$ yields,
\begin{equation}
\doz=\dor\ddor+\frac{M^2r\dor}{R^2}-\frac{M^2r^2}{R^3}+\rho''(d)d\dod
\label{eq:3.z1}
\end{equation}
where we used $(f(d))^\cdot=\rho''(d)d\dod.$ The strong version of the ODE \eqref{eq:3.15} is given by
\begin{equation}
\frac{M^2r}{R}=\dor+R\ddor+M\rho''(d)\dod r,
\label{eq:3.z2}
\end{equation}
which, when multiplied by $\frac{\dor}{R}$ leads to 
\begin{equation}
\rho''(d)d\dod=\frac{M^2r\dor}{R^2}-\frac{\dor^2}{R}-\dor\ddor.
\label{eq:3.z3}
\end{equation}
Substituting, $\rho''(d)d\dod$ in \eqref{eq:3.z1} via \eqref{eq:3.z3} yields,
\begin{equation*}
\doz=-\frac{M^2r^2}{R^3}+\frac{2M^2r\dor}{R^2}-\frac{\dor^2}{R}.
\label{eq:3.z4}
\end{equation*}
The second derivative of $z$ is then given by
\[
\ddoz=-\frac{2M^2r\dor}{R^3}+\frac{3M^2r^2}{R^4}+\frac{2M^2\dor^2}{R^2}+\frac{2M^2r\ddor}{R^2}-\frac{4M^2r\dor}{R^3}-\frac{2\dor\ddor}{R}+\frac{\dor^2}{R^2}.
\label{eq:3.z5}
\]
and the Laplacian becomes
\[
\D z=-\frac{2M^2r^2}{R^4}+\frac{2M\dor^2}{R^4}+\frac{2Mr\ddor}{R^2}-\frac{4Mr\dor}{R^3}-\frac{2\dor\ddor}{R}.
\label{eq:3.z6}
\]
The equations (\ref{eq:3.z3}) and (\ref{eq:3.d2}) allow us to replace the terms with a second derivative $\ddor$
\begin{equation}
\D z=\frac{2}{R^4}[M^2r^2(1+M)+R^2\dor^2(1+M)-2MRr\dor]-\rho''(d)\dod\left(\frac{2M^2r^2}{R^3}-\frac{2d}{R}\right).
\label{eq:3.z7}
\end{equation}
Now define
\begin{equation}
s(R):=M^2r^2(1+M)+R^2\dor^2(1+M)-2RMr\dor.
\label{eq:3.z8}
\end{equation}
Completion of the square, yields
 \[
s(R)=(Mr-R\dor)^2 +M^3r^2+MR^2\dor^2\ge0.
\label{eq:3.z9}
\]
To deal with the `$\rho''$-terms' of (\ref{eq:3.z7}) we  use the form (\ref{eq:3.d3}) of $\dod$ to obtain

\begin{align*}
\left(\frac{2M^2r^2}{R^3}-\frac{2d}{R}\right)\dod=&\left(\frac{2M^2r^2}{R^3}-\frac{2d}{R}\right)\frac{M[(R\dor-r)^2+(M^2-1)r^2]}{R^3+M^2\rho''(d)r^2R}&\\
=&\frac{2M^2[M^3r^4-(2M+M^2)Rr^3\dor+(2+M)R^2r^2\dor^2-R^3r\dor^3]}{R^3(R^3+M^2\rho''(d)r^2R)}&
\label{eq:3.z10}
\end{align*}
Denoting the expression in the brackets by $t(R),$ i.e.
\[
t(R):=M^3r^4-(2M+M^2)Rr^3\dor+(2+M)R^2r^2\dor^2-R^3r\dor^3
\label{eq:3.z11}
\]
we get 
 \begin{eqnarray*}
\D z&=&\frac{2s}{R^4}-\rho''(d)\frac{2M^2t}{R^3(R^3+M^2\rho''(d)r^2R)}\\
&=&\frac{2}{R^4}\frac{[R^3s+\rho''(d)M^2R(r^2s-t)]}{R^3+M^2\rho''(d)r^2R}.
\label{eq:3.z12}
\end{eqnarray*}
Since $\rho''(d)\ge0$ the denominator is nonnegative and we know that $s\ge0,$ therefore, the first term is nonnegative. To complete the proof it is enough to show that
 \begin{eqnarray}
\frac{2}{R^4}\frac{M^2R(r^2s-t)}{R^3+M^2\rho''(d)r^2R}+c_M\doz^2\ge0.
\label{eq:3.z13}
\end{eqnarray}
In order to prove (\ref{eq:3.z13}) we show the slightly stronger statement
 \begin{eqnarray}
2M^2(r^2s-t)+c_M\doz^2R^6\ge0.
\label{eq:3.z14}
\end{eqnarray}
Assume that (\ref{eq:3.z14}) holds, adding the nonnegative term $c_M\doz^2R^3(M^2\rho''(d)r^2R)$ to the left hand side and dividing by $R^3(R^3+M^2\rho''(d)r^2R)$ yields
\[
\frac{2M^2(r^2s-t)+c_M\doz^2R^3(R^3+M^2\rho''(d)r^2R)}{R^3(R^3+M^2\rho''(d)r^2R)}\ge0,
\label{eq:3.z15}
\]
which agrees with (\ref{eq:3.z13}). First we compute the quantities 
\[
\doz R^6 =M^4r^4+R^4\dor^4+(4M^4+2M^2)R^2r^2\dor^2-4M^4Rr\dor
\label{eq:3.z16}
\]
and
\[
2M^2(r^2s-t) =2M^4r^4-2M^2 R^2r^2\dor^2+2M^2R^3r\dor^3+2M^4Rr^3\dor.
\label{eq:3.z17}
\]
Then \eqref{eq:3.z14} becomes
\begin{flalign}
&2M^2(r^2s-t)+c_M\doz R^6& \nonumber\\
&=c_MR^4\dor^4+(2+c_M)M^4r^4+(4M^4c_M+2M^2c_M-2M^2) R^2r^2\dor^2 \nonumber\\
&\hspace{0.5cm}+2M^2(1-2c_M)R^3r\dor^3+2M^4(1-2c_M)Rr^3\dor& \nonumber\\
&=\left(\frac{2c_M-1}{2}\right)(M^{5/4}r-M^{1/4}R\dor)^4+\left(c_M-\left(\frac{2c_M-1}{2}\right)M\right)R^4\dor^4& \nonumber\\
&\hspace{0.5cm}+\left((2+c_M)M^4-\left(\frac{2c_M-1}{2}\right)M^5\right)r^4& \nonumber\\
&\hspace{0.5cm}+\left(4M^4c_M+2M^2c_M-2M^2-6\left(\frac{2c_M-1}{2}\right)M^3\right) R^2r^2\dor^2.&
\label{eq:3.z18}
\end{flalign}
In the last step we completed the quartic form. Note that $\left(\frac{2c_M-1}{2}\right)=\frac{c_M}{M}=\frac{1}{2(M-1)}\ge0.$ Furthermore, 
the coefficients satisfy
\begin{align*}
c_M-\left(\frac{2c_M-1}{2}\right)M=&c_M-\left(\frac{c_M}{M}\right)M=0,&\\
(2+c_M)M^4-\left(\frac{2c_M-1}{2}\right)M^5=&M^4\left((2+c_M)-\left(\frac{c_M}{M}\right)M\right)=2M^4\ge0.&
\end{align*}
The last coefficient in \eqref{eq:3.z18} is nonnegative if it satisfies the following condition
\begin{align*}
4M^4c_M+2M^2c_M-2M^2-6\left(\frac{2c_M-1}{2}\right)M^3=&4M^4c_M+2M^2c_M-2M^2-6c_MM^2&\\
=&4M^4c_M-4M^2c_M-2M^2\ge 0&
\end{align*}
equivalently,
\begin{align*}
1\le(2M^2-2)c_M=M\left(\frac{2M^2-2}{2M-2}\right),
\end{align*}
which is true for all $M\in\N,$ $M\ge2.$ This finishes the proof since all terms in (\ref{eq:3.z18}) are nonnegative yielding $2M^2(r^2s-t)+c_M\doz R^6\ge0$ and with the discussion above  $\D z+(\frac{M}{2(M-1)})\rho''(d)\doz\ge0$ in $(0,1]$.\end{proof}\vspace{0.5cm}

The latter statement now allows to control the behaviour close to the origin. Indeed, $z$ is monotonic close to $0.$ This gives constraints on the quantity $\frac{R\dor}{r},$ which will be useful later.

\begin{lem}\label{lem:3.7} Let $M\in\N,$ $M\ge2.$ and $r$ solves the BVP \eqref{eq:3.0.1}. \\

Then there exists $\d>0$ s.t. $z$ is monotone on $(0,\d).$ Assume, additionally, $r>0$ in $(0,\d).$ Then one of the following conditions holds identically in $(0,\d):$
\begin{align}
\dot{z}\ge0 \;&\mbox{equivalently}\; 0<M^2-M\sqrt{M^2-1}\le\frac{R\dor}{r}\le M^2+M\sqrt{M^2-1},\label{eq:3.23a}&\\
\dot{z}\le0 \;&\mbox{and}\; \left(\frac{r}{R}\right)^\cdot> 0 \;\mbox{or}\label{eq:3.23b}&\\
\dot{z}\le0 \;&\mbox{and}\; \left(\frac{r}{R}\right)^\cdot< 0.&
\label{eq:3.23c}
\end{align} 
\label{lem:3.z}
\end{lem}
\begin{proof}
By lemma \ref{lem:3.4} we know either that $z(x)$ is constant, in which case the monotonicity is given, or $z$ does not attain a maximum in $B\sm\{0\}.$ It follows that $R\mapsto z(R)$ can only have one local minimum in $(0,1],$ and therefore $\doz$ can only change sign once. Hence, there exists a $\d>0$ s.t. $R\mapsto z(R)$ is monotone on $(0,\d).$\\   

Suppose now that $\d>0$ s.t.\! the above holds and $r(R)>0$ for all $R\in(0,\d)$ and recall $z=\frac{\dor^2}{2}+\frac{M^2r^2}{2R^2}+f(d).$ Then the derivative of $z$ is given by
\begin{flalign*}
\dot{z}=\dor\ddor+\frac{M^2r\dor}{R^2}-\frac{M^2r^2}{R^3}+(f(d))^\cdot
\label{eq:3.23d}
\end{flalign*} 
Using the definition of $f$ we obtain
\[
(f(d))^\cdot=(\rho'(d)d-\rho(d))^\cdot=\rho'(d)\dod+\rho''(d)d\dod-\rho'(d)\dod=\rho''(d)d\dod.
\label{eq:3.23e}
\]
On the other hand, since $r$ is a strong solution to the ODE \eqref{eq:3.15} in $(0,1],$ all derivatives exist in a strong sense. In particular,
\begin{align*}
\frac{M^2r}{R}+M\rho'(d)\dot{r}=&\left(R\dot{r}+M\rho'(d)r\right)^\cdot& \nonumber\\
=&\dor+R\ddor+M\rho'(d)\dot{r}+M\rho''(d)\dod r.&
\label{eq:3.23f}
\end{align*}
Rearranging the above equations and multiplying it by $\frac{\dor}{R}$ yields
\[
(f(d))^\cdot=\rho''(d)d\dod=\frac{M^2r\dor}{R^2}-\frac{\dor^2}{R}-\dor\ddor.
\label{eq:3.23g}
\]
Therefore,
\begin{align}
\dot{z}(R)=-\frac{r^2}{R^3}\left[\left(\frac{R\dor}{r}\right)^2-2M^2\left(\frac{R\dor}{r}\right)+M^2\right].
\label{3.23h}
\end{align}

This polynomial of variable $\frac{R\dor}{r}$ has roots at $\lambda_{\pm}=M^2\pm M\sqrt{M^2-1}.$ It is negative between these roots and positive otherwise. In particular, if $\doz\ge0$ then , from \eqref{3.23h}, it follows that $M^2-M\sqrt{M^2-1}\le\frac{R\dor}{r}\le M^2+M\sqrt{M^2-1},$ which agrees with \eqref{eq:3.23a}. Now for \eqref{eq:3.23b} and \eqref{eq:3.23c} from $\dot{z}\le0$ on $(0,\d)$ it follows $\frac{R\dor}{r}\le M^2- M\sqrt{M^2-1}<1$ or $\frac{R\dor}{r}\ge M^2+ M\sqrt{M^2-1}>1.$ Then the calculation
\begin{eqnarray*}
\frac{R\dor}{r}&\gtrless& 1 \nonumber\\
\Rightarrow\frac{\dor}{R}&\gtrless& \frac{r}{R^2} \nonumber\\
\Rightarrow\left(\frac{r}{R}\right)^\cdot&=&\frac{\dor}{R}-\frac{r}{R^2}\gtrless 0
\label{3.23i}
\end{eqnarray*}
shows (\ref{eq:3.23b}) and (\ref{eq:3.23c}).\end{proof}\vspace{0.5cm}

Up to this point we have narrowed down the number of possibilities, how $d$ and $z$ can behave close to $0,$ enough so that there are only a few cases left, which now can be treated individually. The shape of $d$ and $z$ will demand certain conditions on $r,$ which can either be matched by $r$ or will lead to contradictions, excluding these cases. This will reduce the number of types even further and leaves only the following restrictive statement:   
\begin{lem}Let $M\in\N,$ $M\ge2$ and suppose $r$ solves the BVP \eqref{eq:3.0.1}. Then $d(0)=0.$
Moreover, one of the two situations occurs:\\
i) $r$ is lifting-off delayed, i.e. there exists, $0<\d<1$ s.t. $r\equiv0$ on $[0,\d].$ Then $r\in C^\infty([0,1]).$\\
ii) $r$ is lifting-off immediately, i.e. $r,\dor>0$ away from zero. Then $r\in C^1([0,1])$ and $\dor(0)=0.$
\label{lem:3.2.100}
\end{lem}
\begin{proof}
We show that only $d(0)=0$ is possible and either $r\equiv 0$ near zero or $r>0$ away from zero and $\doz\ge0.$ All other situations are excluded by contraposition.\\
   
Recall, that $\lim\limits_{R\ra0} d(R)$ exists and agrees with $d(0)$ due to continuity, which was proven in lemma \ref{lem:Mdet}. Also notice that lemma \ref{lem:3.7} guarantees that $\lim\limits_{R\ra0} z(R)$ makes sense, however the limit might be $+\infty.$ In particular, we know $0\le\lim\limits_{R\rightarrow0} z(R)\le+\infty$ by the definition of $z$ and the behaviour of $f.$ \\

1. Case: $\lim\limits_{R\rightarrow0} d(R)=:l\in(0,\infty).$\\
By continuity of $d$ on $(0,1]$ there exists $\delta>0$ s.t. $r\dor>0$ on $(0,\delta)$ implying $r,\dor>0$ on $(0,\d)$ (even on $(0,1]$, by monotonicity).\\
First, assume $\lim\limits_{R\rightarrow0} z(R)=+\infty.$ In this case only $\doz\le0$ near zero is possible implying $\frac{r}{R}$  to be monotone on $(0,\d).$ The mean value theorem guarantees the existence of the two sequences $R_j'\ra0$ and $R_j\in(0,R_j')$ for any $j\in \N$ s.t. 
\begin{equation}
\lim\limits_{j\rightarrow\infty} \dor(R_j)=\lim\limits_{j\rightarrow\infty} \frac{r(R_j')}{R_j'}=:m.
\label{eq:3.l.9.1}
\end{equation}
By lemma \ref{lem:3.7} we know that there are two different cases, either $\dor>\frac{r}{R}$ or $\dor<\frac{r}{R}$ on $(0,\d).$ In the first case by $\dor>\frac{r}{R}$ and \eqref{eq:3.l.9.1} we have
\[
l=\lim\limits_{j\ra\infty} d(R_j')=\lim\limits_{j\ra\infty} \frac{M\dor(R_j') r(R_j')}{R_j'}\ge\lim\limits_{j\ra\infty} M \frac{r^2(R_j')}{R_j'^2}=Mm^2.
\]
Hence, we know $m\le\sqrt{\frac{l}{M}}<+\infty$ and together with the property $\dor>\frac{r}{R}$ on $(0,\d)$ it holds
\[
\lim\limits_{j\rightarrow\infty} z(R_j)=\lim\limits_{j\rightarrow\infty}\left[\frac{\dor^2(R_j)}{2}+\frac{M^2r^2(R_j)}{2R_j^2}+f(d(R_j))\right]\le\left[\frac{1}{M}+M\right]\frac{l}{2}+f(d(l))<+\infty,\]
contradicting $\lim\limits_{R\rightarrow0} z(R)=+\infty.$

In the 2nd case, that is $\dor<\frac{r}{R}$ on $(0,\d),$ it holds
\[
l=\lim\limits_{j\ra\infty} d(R_j)=\lim\limits_{j\ra\infty} \frac{M\dor(R_j) r(R_j)}{R_j}\ge\lim\limits_{j\ra\infty} M\dor^2(R_j)=Mm^2.
\]
Again $m\le\sqrt{\frac{l}{M}}<+\infty$ and together with the property $\dor<\frac{r}{R}$ on $(0,\d)$ we have
\[
\lim\limits_{j\rightarrow\infty} z(R_j')=\lim\limits_{j\rightarrow\infty}\left[\frac{\dor^2(R_j')}{2}+\frac{M^2r^2(R_j')}{2R_j'^2}+f(d(R_j'))\right]\le\left[\frac{1}{M}+M\right]\frac{l}{2}+f(d(l))<+\infty,\]
contradicting $\lim\limits_{R\rightarrow0} z(R)=+\infty.$\\
\vspace{0.5cm}

Now assume the limit exists, i.e. $\lim\limits_{R\rightarrow0} z(R)=:n\in[0,\infty).$ For the sake of contradiction we show that the right limit of $\dor$ exists in $0$ and that it is nonzero.\\
Introduce the new variables $\nu_1(R):=\dor(R)$ and $\nu_2(R):=\frac{Mr(R)}{R}.$ Then we can interpret the functions $d$ and $z$ as functions depending on these new variables $d(\nu_1,\nu_2)=\nu_1\nu_2$ and $z(\nu_1,\nu_2)=\frac{\nu_1^2}{2}+\frac{\nu_2^2}{2}+f(\nu_1\nu_2)$ on the set
$\mathcal{V}=\{(\nu_1,\nu_2):\nu_1>0,\nu_2>0\}.$\\
Consider 
\[
K:=\{(\nu_1,\nu_2):\frac{l}{2}\le d\le 2l\}\cap\{(\nu_1,\nu_2):s-1\le z\le s+1\}
\]
here $l>0$ and $s>1$ are parameters, in particular, $s$ is not the function introduced in \eqref{eq:3.z8}. $K$ is a compact subset of $\mathcal{V}$ due to the continuity of $d$ and $z.$ In particular, $\{d=l\}\cap\{z=s\}$ consists of at most two points $\{(a,b),(b,a)\}\ss\mathcal{V}.$\\ Now again by the continuity of $r$ and $\dor$ we know that for all $\ve>0$ there exists $R_0(\ve)>0$ s.t.\@ for all $0<R<R_0,$ $(\nu_1(R),\nu_2(R))\in B_\ve(a,b)\cup B_\ve(b,a)\ss K.$ Now we need to show that $R\mapsto(\nu_1(R),\nu_2(R))$ remains in one of these balls for all $0<R<R_0.$ If $a=b$ this is immediately true. So suppose $a\not=b.$ Then we can choose $\ve>0$ so small that the balls become disjoint, i.e. $B_\ve(a,b)\cup B_\ve(b,a)=\emptyset.$ Recall that $r\in C^\infty((0,1])$ therefore, the curve $R\mapsto(\nu_1(R),\nu_2(R))=(\dor(R),\frac{Mr(R)}{R})$ is connected and remains in one ball, say $B_\ve(a,b).$ Since $\ve>0$ was arbitrary, we see that
\begin{eqnarray*}
\lim\limits_{R\rightarrow0}\dor(R)=a\in(0,+\infty)
\end{eqnarray*}
Hence, $r\in C^1([0,1])$ and $\dor(0)>0.$\\

Consider the rescaled function $r_\ve(R):=\ve^{-1}r(\ve R)$ for $0<\ve<1$ with the derivatives $\dor_\ve(R)=\dor(\ve R)$ and $\ddor_\ve(R)=\ve\ddor(\ve R)$ for all $R\in(0,1].$ Since $r$ solves the ODE (\ref{eq:3.15}) strongly in $(0,1)$ so does the rescaled version. Indeed note that 
\begin{align*}
d_{r_\ve}(R)=&\frac{Mr_\ve(R)\dor_\ve(R)}{R}=\frac{Mr(\ve R)\dor(\ve R)}{\ve R}=d_r(\ve R)\;\;\mbox{and}&\\
\dod_{r_\ve}(R)=&(d_r(\ve R))^\cdot=\ve \dod_r(\ve R).&
\end{align*}
Hence,
\begin{eqnarray*}
\dor_\ve(R)+R\ddor_\ve(R)+M\rho''(d_{r_\ve}(R))\dod_{r_\ve}(R)r_\ve(R)-\frac{M^2r_\ve(R)}{R}=\nonumber\\
=\dor(\ve R)+(\ve R)\ddor(\ve R)+M\rho''(d_r(\ve R))\ve \dod_r(\ve R)\ve^{-1}r(\ve R)-\frac{M^2r(\ve R)}{\ve R}=0.
\end{eqnarray*}
where the last equality holds since it agrees with the strong form of the ODE (\ref{eq:3.15}) evaluated at $\ve R.$\\
Now if $\ve\ra 0,$ then the rescaled function $r_\ve$ converges uniformly to the linear map $r_0(R):=aR$, i.e. $r_\ve\ra aR$ and $\dor_\ve\ra a$ uniformly in $[0,1].$ But then the weak form of the ODE of $r_\ve$ converges to the weak form of $r_0:$
\begin{eqnarray*}
0&=&\lim\limits_{\ve\ra0}\int\limits_{0}^1{\dot{g}\left(\dor_\ve R+M\rho'(d_{r_\ve})r_\ve\right)+g\left(\frac{M^2r_\ve}{R}+M\rho'(d_{r_\ve})\dor_\ve\right)\;dR}\nonumber\\
&=&\int\limits_{0}^1{\dot{g}\left(\dor_0(R)R+M\rho'(d_{r_0})r_0\right)+g\left(\frac{M^2r_0}{R}+M\rho'(d_{r_0})\dor_0\right)\;dR}.
\end{eqnarray*}
Therefore, $r_0$ is a weak solution to the ODE (\ref{eq:3.15}). Since $r_0$ is smooth it also needs to satisfy (\ref{eq:3.z2}), which is not true since plugging $r_0$ into (\ref{eq:3.z2}) yields
\[M^2a=a,\]
which is not satisfied since $M\ge2,\;a>0.$\\
 
2. Case: $\lim\limits_{R\rightarrow0} d(R)=0.$ \\
There are only two possible scenarios: Either $r\equiv 0$ in $[0,\d]$ or $r>0$ in $(0,\d).$

If $r\equiv0$ near zero then $\dor\equiv0$ in $(0,\d]$ and we can easily see $\dor(0)=0$ and $\dor\in C^1([0,1]).$ Moreover, this argument works for all derivatives of $r$ yielding $r\in C^\infty([0,1]).$\\

Assume instead that $r$ lifts-off immediately, i.e. there exists $\d>0$ s.t. $r>0$ on $(0,\d).$ Then lemma \ref{lem:3.z} holds, assume first $\doz\le0.$ Again, $\frac{r}{R}$ is monotone on $(0,\d)$ and the mean value theorem implies the existence of of the two sequences $R_j'\rightarrow0$ and $R_j\in(0,R_j')$ for any $j\in \N$ s.t. 
\[
\lim\limits_{j\rightarrow\infty} \dor(R_j)=\lim\limits_{j\rightarrow\infty} \frac{r(R_j')}{R_j'}=:m.
\]
Assume now $\dor>\frac{r}{R}$ on $(0,\d).$ Then \[0=\lim\limits_{R\rightarrow0} d(R)=\lim\limits_{j\rightarrow\infty}\frac{M\dor(R_j')r(R_j')}{R_j'}\ge\lim\limits_{j\rightarrow\infty}\frac{Mr^2(R_j')}{R_j'^2}=Mm^2.\] 
Hence, $m=0.$ Now, by $\dor>\frac{r}{R}$ on $(0,\d)$ we get $\dor(R_j)\ra0$ and $\frac{r(R_j)}{R_j}\ra0$ and therefore $\lim\limits_{R\rightarrow0} d(R_j)=0$ if $j\ra\infty.$ This yields
\[
\lim\limits_{j\rightarrow\infty} z(R_j)=\lim\limits_{j\rightarrow\infty}\left[\frac{\dor^2(R_j)}{2}+\frac{M^2r^2(R_j)}{2R_j^2}+f(d(R_j))\right]=0.\]
Since $z\ge0$ and by assumption $\doz\le0$ on $[0,\d)$ we have $z=0$ on $[0,\d).$ By the non-negativity of $d\ge0$ and $f(d)\ge0$ (following from $f(0)=0$ and recalling $(f(d))^\cd=\rho''(d)d\dod\ge0$) we finally know that $r=\dor=0$ on $[0,\d),$ 
contradicting the assumption that $r$ lifts off immediately. One can argue similarly in the case when $\dor<\frac{r}{R}$ on $(0,\d).$ \\

Finally, assume $\doz\ge0.$ Then $(\ref{eq:3.23a})$ holds, i.e. $\dor\sim\frac{r}{R}$ on $(0,\d).$ Since $d(R)\rightarrow0$ for $R\rightarrow0$ 
\[
\dor(0)=\lim\limits_{R\rightarrow0} \dor(R)=\lim\limits_{R\rightarrow0}\frac{r(R)}{R}=0
\label{eq:}
\]
and again $r\in C^1([0,1])$ with $\dor(0)=0.$\end{proof}

We end this paragraph by showing that the constructed maps $u=re_{MR}$ s.t. $r$ solves the BVP \eqref{eq:3.0.1} are stationary points of the functional \eqref{eq:1.1}.

\begin{lem} Let $u\in\A_r^M$ with $u=r e_{M\th}$ s.t. $r$ solves the BVP \eqref{eq:3.0.1}. Then $u$ solves the ELE \eqref{eq:ELE1.1} weakly, i.e. 
\[\int\limits_B{\grad_\xi W(\grad u)\cd\grad \vp \;dx}=0 \;\mb{for all}\; \vp\in C_c^\infty(B,\R^2).\]
\label{lem:2.7}
\end{lem}
\begin{proof}
Lemma 3.6 of \cite{BOP91MS} applies and shows that $u$ solves the ELE strongly in $B\sm\{0\},$ this can be reformulated in the following sense, $u$ satisfies
\[\int\limits_B{\grad_\xi W(\grad u)\cd\grad \vp \;dx}=0\]
for all $\vp\in C_c^\infty(B,\R^2)$ with $\vp\equiv0$ near the origin.\\ Now we can follow the strategy of Theorem 3.11 in \cite{BOP91MS} to upgrade this to arbitrary test functions. For this sake, take $\eta_\ve\in C^\infty(B)$ s.t. $\eta_\ve\equiv0$ on $B_\ve$ and $\eta_\ve\equiv1$ on $B\sm B_{2\ve},$ $0\le \eta_\ve \le1$ and there exists $c>0$ s.t. $|\grad \eta_\ve|\le \frac{c}{\ve}.$ Take an arbitrary test function $\psi\in C_c^\infty(B,\R^2)$ and set $\vp=\eta_\ve\psi$ then $\vp$ vanishes close to the origin. Hence,
\begin{align*}
0&=\int\limits_B{\grad_\xi W(\grad u)\cd\grad \vp \;dx}&\\
&=\int\limits_B{\grad_\xi W(\grad u)\cd\eta_\ve\grad \psi \;dx}+\int\limits_B{\grad_\xi W(\grad u)\cd(\grad\eta_\ve\ot \psi) \;dx}&
\end{align*}

Then the first integral converges:
\[
\lim\limits_{\ve\ra0}\int\limits_B{\grad_\xi W(\grad u)\cd\eta_\ve\grad \psi \;dx}=\int\limits_B{\grad_\xi W(\grad u)\cd\grad \psi \;dx}
\]
by dominated convergence. 
Since, $\grad u\in C^0$ there exists $C>0$ s.t. $\|\grad_\xi W(\grad u)\|_{C^0}\le C.$ Then
\[
\left|\int\limits_B{\grad_\xi W(\grad u)\cd(\grad\eta_\ve\ot \psi) \;dx}\right|\le C\|\psi\|_{C^0}\int\limits_B{|\grad\eta_\ve| \;dx}\le \frac{C}{\ve}\|\psi\|_{C^0}\La^2(B_{2\ve})
\]
Hence, the second integral vanishes for $\ve\ra0$ and 
\[\int\limits_B{\grad_\xi W(\grad u)\cd\grad \psi \;dx}=0\]
holds for every $\psi\in C_c^\infty(B,\R^2).$\end{proof}

\section{Advanced BOP-theory}
\label{sec:3.3}
In the following, we want to investigate if these stationary points are even more regular. A first step in that direction is the next lemma. This has been observed for the BOP-case by Yan and Bevan, see \cite[Lem 3.(i)]{BY07}  and \cite[Lem 1.(i)]{Y06}. 
\begin{lem} Let $M\in\N,$ $M\ge1$ and suppose $r$ solves the BVP \eqref{eq:3.0.1}. Then\\
(i) $\liminf\limits_{R\ra0}\ddor(R)\ge0$ and\\
(ii) $\lim_{R\ra0}\ddor(R)R=0.$
\label{lem:3.1}
\end{lem}
\begin{proof}
(i) Assume not, then $\liminf\limits_{R\ra0}\ddor(R)<0.$ But then there exits a small interval s.t. $\ddor<0$ on $(0,\d).$ By the mean value theorem it follows the existence of some $\xi\in(0,\d)$ , s.t.
 \[0>\ddor(\xi)=\frac{\dor(\d)}{\d}\ge0,\]
which is a contradiction.\\ 
(ii) Recall the ODE
\[
\ddor(R)R=\frac{M^2r}{R}-\dor-M\rho''(d)\dod r
\]
Then 
\[
0\le\liminf\limits_{R\ra0}{\ddor(R)R}\le\limsup\limits_{R\ra0}{\ddor(R)R}\le\limsup\limits_{R\ra0}{\frac{M^2r(R)}{R}}=0.
\]
This shows the second claim. \end{proof}\vspace{1cm}

\subsection{Delayed lift-off solution for arbitrary $\rho$}
From now on, we will distinguish between the different shapes of the solutions to the BVP \eqref{eq:3.0.1}. For this recall, that $r$ is an immediate lift-off solution if $r,\dor>0$ away from zero.  We will denote such solutions by $r_0,$ while we will call delayed lift-off solutions by $r_\d$ for $0<\d<1.$ As a reminder, $r_\d$ is a delayed lift-off solution if there is $0<\d<1$ s.t. $r\equiv0$ on $[0,\d].$ The $\d$ indicates the point such that $r_\d\equiv0$ in $[0,\d]$ but also $r_\d(R)>0$ for all $R\in(\d,1].$ The first statement will be that an immediate lift-off solution $r_\d$ is zero up to $\d$ and that it needs to solve the BVP \eqref{eq:3.3.3} below. This fact will be crucial for the uniqueness result in lemma \ref{Lem:3.3.3}.
\begin{lem}Let $0<\ga<\infty.$ 
If for some $\d>0,$ $r_\d$ solves the BVP \eqref{eq:3.0.1} then \[r_{\d}(R)=\left\{\begin{array}{ccc}
0& {\mbox{in}}& (0,\d],\\
\tilde{r}_{\d,0}& {\mbox{in}}& (\d,1],
\end{array}
\right.\label{eq:3.3.2}\]
where $\tilde{r}_{\d,0}\in C^\infty((\d,1))\cap C^0([\d,1])$ is the unique solution of
\begin{align}\left\{\begin{array}{ccc}
Lr=M\rho''(d)\dod r & {\mbox{in}}& (\d,1),\\
r(\d)=0 , &r(1)=1.&
\end{array}
\right.\label{eq:3.3.3}\end{align}
\end{lem}
\begin{proof}
Assume $r_\d$ is a solution to the BVP \eqref{eq:3.0.1}. Then, $r_\d=0$ on $[0,\d]$ by definition and $r_\d|_{[\d,1]}\in C^{\infty}([\d,1])$ is the unique solution to 
\[\left\{\begin{array}{ccc}
Lr=M\rho''(d)\dod r & {\mbox{in}}& (\d,1),\\
r^{(k)}(\d)=0\;\mb{for all}\;k\in\N , &r(1)=1&
\end{array}
\right.\]
Hence, $r_\d|_{[\d,1]}$ solves \eqref{eq:3.3.3}.\end{proof}\vspace{0.5cm}

Furthermore, we demonstrate that such a solution, if it exists, needs to be unique. 
\begin{lem}Let $0<\ga<\infty.$ Assume there exists a solution to the BVP \eqref{eq:3.0.1} of the form $r_\d$, $\d>0.$ Then there exists a unique $\d>0$ in the sense that
\[r_{\d}(R)=\left\{\begin{array}{ccc}
0& {\mbox{in}}& (0,\d],\\
\tilde{r}_{\d,0}& {\mbox{in}}& (\d,1]
\end{array}
\right.\]
and $\tilde{r}_{\d,0}$ lifts off immediately, i.e. $\tilde{r}_{\d,0}(R)>0$ for all $R\in(\d,1].$
\label{Lem:3.3.3}
\end{lem}
\begin{proof}
We can always choose $\d$ to be maximal in the sense that $r_\d\equiv0$ on $[0,\d]$ and $\tilde{r}_{\d,0}$ lifts off immediately. For the uniqueness assume that $r_\d$ and $r_{\d'}$ are two solutions to \eqref{eq:3.0.1} for $0<\d\le\d'<1$. Then both need to satisfy  \[\left\{\begin{array}{ccc}
Lr=M\rho''(d)\dod r & {\mbox{in}}& (\d,1),\\
r(\d)=0\;, &r(1)=1.&
\end{array}
\right.\]
The Picard-Lindelöf Theorem\footnote{see, \cite{GT12} Theorem 2.5, Corollary 2.6.} yields, $r_\d\equiv r_{\d'}$ on $[\d,1]$, trivially on the complete interval $[0,1].$ Hence, $\d=\d'$ by maximality.\end{proof}

\textbf{Conclusion:}
For arbitrary behaviour of $\rho$ and $0<\ga<\infty$ the above discussion shows that (smooth) delayed lift-off solutions to the BVP \eqref{eq:3.0.1} either do not exist or there is at most one. The latter lemma shows that the representation and $\d$ may vary, but the solution remains the same.

\subsection{Delayed lift-off $\rho$}

Lets assume $\rho$ itself is a delayed lift-off function, i.e. there exists $\tilde{s}>0$ s.t.  $\rho\equiv0$ on $[0,\tilde{s}]$. Since, $d(0)=0$ and $d$ is continuous there exists a small interval $[0,\d)$ s.t. $R\mt\rho(d(R))\equiv0.$ So the $\rho-$term of the considered functional \eqref{eq:1.1} vanishes reducing the functional to the Dirichlet energy, at least close to the origin. It is then well known, that stationary points to the Dirichlet energy are harmonic functions and therefore smooth. Giving hope, that in case of a delayed lift-off $\rho,$ all solutions are indeed smooth. The following result states that this is exactly true. Notice, that it will be enough to consider immediate lift-off solutions $r_0$ since, by the previous discussion, we already know that delayed lift-off solutions have to be smooth.
     
\begin{lem}Let $0<\ga<\infty$ and assume there exists $\tilde{s}>0$ s.t. $\rho(s)=0$ for all $s\in[0,\tilde{s}]$ and $\rho(s)>0$ for all $s>\tilde{s}.$ Further suppose there exists a solution of the BVP \eqref{eq:3.0.1} of the form $r_0.$\\
If $r_0$ solves \eqref{eq:3.0.1} then there exists a unique $\d=\d(\tilde{s},r_0,\dor_0)>0$ and some $0<a<1$ s.t. 
\begin{align}r_{0}(R)=\left\{\begin{array}{ccc}
a\left(\frac{R}{\d}\right)^M& {\mbox{in}}& (0,\d],\\
\tilde{r}_{\d,a}& {\mbox{in}}& (\d,1]
\end{array}\label{eq:3.3.4}
\right.\end{align}
and $\tilde{r}_{\d,a}$ the unique smooth solution of
\begin{align}\left\{\begin{array}{ccc}
Lr=M\rho''(d)\dod r & {\mbox{in}}& (\d,1),\\
r(\d)=a , &r(1)=1.&
\label{eq:3.3.5}
\end{array}
\right.\end{align}
and $\tilde{r}_{\d,a}$ is not in the kernel of $L$ for at least a short period of time, i.e. there exists $\ve>0$ s.t. $L\tilde{r}_{\d,a}(R)>0$ for all $R\in(\d,\d+\ve].$ Moreover, $(a\left(\frac{R}{\d}\right)^M)^{(k)}(\d)=\tilde{r}_{\d,a}^{(k)}(\d)$ for all $k\in \N.$
\label{lem:3.4.1}
\end{lem}
\begin{proof}
Since $d(0)=0$ and $d\in C([0,1])$ and the delay of $\rho$ there exists a unique $\d=\d(\tilde{s},r_0,\dor_0)>0$ s.t. $R\mapsto\rho(d(R))\equiv0$ for all $R\in[0,\d]$ and $\rho(d(R))>0$ for all $R>\d.$ Then $d(R)>0$ for $\d<R\le1$ and there exists an $\ve>0$ s.t. $\dod(R)>0$ for $R\in(\d,\d+\ve].$ Hence, $L\tilde{r}_{\d,a}(R)=M\rho''(d)\dod r>0$ for $R\in(\d,\d+\ve].$
Then $r_0$ needs to solve the following ODE
\[\left\{\begin{array}{ccc}
Lr=0& {\mbox{in}}& (0,\d),\\
r(0)=0,&r(\d)=a&
\end{array}
\right.\]
for some $0<a<1.$ Indeed $a$ can not exceed $1$ $(a>1)$ since $\dor\ge0$ and $r(1)=1.$ If $a=1$ then $r\equiv1$ in $[\d,1]$ implying $d\equiv0$ in $[0,1],$ a contradiction. $a=0$ is excluded by the assumption that the solution $r_0$ is an immediate lift-off function.\\
Then $r_0|_{[\d,1]}$ solves \eqref{eq:3.3.5} uniquely and $r_0$ takes the form \eqref{eq:3.3.4}.
Since $r_0\in C^\infty((0,1])$ all derivatives of $r_0$ need to agree at $\d,$ i.e. $(a\left(\frac{R}{\d}\right)^M)^{(k)}(\d)=\tilde{r}_{\d,a}^{(k)}(\d)$ for all $k\in \N.$\end{proof}

\begin{re} By construction $\d$ depends not only on $\tilde{s}$ but also on $r_0,\dor_0$. For every $r_0$ the $\d$ may vary, destroying any chance for a uniqueness result similar to Lemma \ref{Lem:3.3.3}.
\end{re}

\textbf{Conclusion:}
We can not guarantee the existence of solutions to the ODE of the form $r_0.$ But if they exist, they have to be smooth.
Moreover, combining this with our knowledge of delayed-lift off solutions, we are able to conclude that in case of a delayed $\rho$ all stationary points in the class $\A_r^M$ are smooth.\\

\subsection{Immediate lift-off $\rho$}

If $\rho$ is an immediate lift-off function ($\rho(s)>0$ for all $s>0$), then we are not (yet) able to show, that $r_0$ solutions need to be any smoother then $C^1.$ However, as a next statement we show that under the additional assumption $u=r_0e_{MR}\in W^{2,2}(\ol{B},\R^2)$ full smoothness of $r,d$, and $u$ can be obtained. Additionally, we then describe a necessary condition \eqref{eq:3.3.200}, which needs to be satisfied for any smooth enough $r_0.$ As a consequence, the limit of the quantity $\frac{R\dor_0}{r_0}$ if $R$ tends to $0,$ playing a key role in lemma \ref{lem:3.7},  can then finally be determined.\\

\begin{lem} (Higher-order regularity)\label{lem:3.3.200} Let $0<\ga<\infty$ and $\rho(d)>0$ for all $d>0.$ Assume there exists a solution to the BVP \eqref{eq:3.0.1} of the form $r_0$ and assume $u=r_0e_{MR}\in W^{2,2}(\ol{B},\R^2).$\\

Then $r_0\in C^{\infty}([0,1]),$ the corresponding $d\in C^{\infty}([0,1]),$ and $u=r_0e_{MR}\in C^{\infty}(\ol{B},\R^2).$ In particular, it holds that $r_0\in C^{1,\al}([0,1])$ for any $\al\in(0,1)$ and there exists $\d>0$ s.t.
\begin{equation}M^2-C_{\al,M}\rho''(d_{r_0})R^{2\al}\le\frac{R\dor_0}{r_0}+\frac{R^2\ddor_0}{r_0}< M^2 \;\mb{for all}\; R\in (0,\d)
\label{eq:3.3.200}\end{equation}
and as a consequence we have \[D_M:=\lim\limits_{R\ra0}\frac{R\dor_0}{r_0}=M.\] 
\end{lem}
\begin{proof}
In the following we will suppress $r_0$ in favour of $r.$\\ 

By the higher-order regularity result \cite[Thm 1.4]{D2} we know that $u=re_{MR}\in W^{2,2}(\ol{B},\R^2)$ and $d\in C^{1}([0,1]),$ is enough to imply full local smoothness $u=re_{MR}\in C_{loc}^{\infty}(B,\R^2)$ which guarantees smoothness at the origin, the smoothness can then be extended by limit taking of $r,$ $d$ and their derivatives up to the harmless boundary at $R=1,$ showing the claim.\\

Now for the second part of the statement we first note that $\rho(d)>0$ for all $d>0$ implies that there exists a $\d>0$ s.t. $R\mapsto\rho''(d(R))>0$ for all  $R\in(0,\d).$ This yields
\begin{equation}0<M\rho''(d)=\frac{Lr}{\dod r},\;\mb{on}\;(0,\d). \label{eq:3.3.98}\end{equation} 
Since, $r,\dor>0$ in $(0,1]$ we can infer $\dod>0$ and $Lr>0$ on $(0,\d).$
Since, $r\in C^{1,\al}([0,1])$ there exists $c_\al>0$ and $\d>0$ s.t.  
\begin{equation}|\ddor(R)|\le c_\al R^{\al-1},\;\mb{for all}\;R\in(0,\d). \label{eq:3.3.99}\end{equation}
Assume not. Then for all $\d>0,$ $c>0$ there exists $R\in (0,\d)$ s.t.
\begin{equation}|\ddor(R)|> c R^{\al-1}. \label{eq:3.3.100}\end{equation}
Fix $\d>0.$ Then for all $c>0$ take $R_c\in (0,\d)$ s.t. the latter inequality holds. By continuity of $\ddor$ in $(0,1],$ there exists an $\ve>0$ s.t. \eqref{eq:3.3.100} holds even for all $R'\in(R_c-\ve,R_c+\ve).$ By integration we get
\[|\dor(R_c+\ve)-\dor(R_c)|=\int\limits_{R_c}^{R_c+\ve}{|\ddor(R')|\;dR'}>\frac{c}{\al}((R_c+\ve)^\al-R_c^\al).\]
Hence, for all $c>0$ we can find $R_c\in (0,\d)$ s.t. $\dor$ is not Hölder continuous at $R_c$ with Hölder constant $\frac{c}{\al}.$ Since $c>0$ is arbitrary this contradicts $r\in C^{1,\al}([0,1]).$\\ 

Then \eqref{eq:3.3.99} implies
\[0<\dod r=M\left(\frac{r\dor^2}{R}+\frac{r^2\ddor}{R}-\frac{r^2\dor}{R^2}\right)\le M\frac{r}{R}\left(\dor^2+r|\ddor|\right)\le c_{\al}^2\al^{-2}M\frac{r}{R}R^{2\al}\]
for all $R\in(0,\d).$ Together with \eqref{eq:3.3.98} we get
\[0<\frac{R}{r}Lr\le C_{\al,M}\rho''(d)R^{2(k-1+\al)} \;\mb{on}\;(0,\d),\]
where $C_{\al,M}:=M^2c_{\al}^2\al^{-2}>0.$
Using the explicit form of $Lr>0$ yields the claimed inequalities
\[M^2-C_{\al,M}\rho''(d)R^{2\al}\le\frac{R\dor}{r}+\frac{R^2\ddor}{r}< M^2 \;\mb{near}\; 0.\]
Taking the limit $R\ra0$ yields 
\[\lim\limits_{R\ra0}\left(\frac{R\dor}{r}+\frac{R^2\ddor}{r}\right)=M^2.\]
With the notation $D_M:=\lim\limits_{R\ra0}\frac{R\dor}{r}$ and $E_M:=\lim\limits_{R\ra0}\frac{R\ddor}{\dor}$ we get by L'H\^{o}pital's rule
\[D_M=\lim\limits_{R\ra0}\frac{R\dor}{r}=\lim\limits_{R\ra0}\frac{R\ddor+\dor}{\dor}=1+E_M\] and
\begin{align*}M^2&=\lim\limits_{R\ra0}\left(\frac{R\dor}{r}+\frac{R^2\ddor}{r}\right)=D_M+\lim\limits_{R\ra0}\left(\frac{R\dor}{r}\frac{R\ddor}{\dor}\right)&\\
&=D_M+D_ME_M=D_M+D_M(D_M-1)=D_M^2,&
\end{align*}
hence $D_M=M.$
\end{proof}

\begin{re}
The necessary condition given in \eqref{eq:3.3.200} is relevant to our discussion, since Improving the lower bound, in  up to the point where it matches the upper one, would show that $r$ can't be of class $C^{1,\al}.$ However, we don't know if it is indeed true or how to show it.
\end{re}

\textbf{Conclusion:}
In this case there are two possibilities: There is at most one smooth delayed lift-off solution $r_\d$ for some $\d>0$ or there could be an immediate $C^1-$lift-off solution of the form $r_0,$ remaining open if the regularity of $r_0$ must be any better.\\

\textbf{Proof of theorem \ref{thm:1.3 MT}:}
\begin{enumerate}
\item See lemma \ref{lem:2.7}.
\item The $M=1-$case is worth mentioning, however fairly standard and amounts to the statement that the integrand is strictly quasiconvex.
\begin{enumerate}
\item Consequence of (b).
\item Follows from the lemmatas \ref{lem:2.1}, \ref{lem:2.3}, and \ref{lem:3.2.100}.
\item See the lemmatas \ref{lem:Mdet} and \ref{lem:3.2.100}.
\item See the lemmatas \ref{lem:3.2.100} and \ref{lem:3.4.1}.
\item See the lemmatas \ref{lem:3.2.100} and \ref{lem:3.3.200}, completing the proof of the theorem and the paper.
\end{enumerate}
\end{enumerate}
\qed

\bibliography{LiteraturePhDMD}

\begin{thebibliography}{10}

\bibitem{AF87}
Emilio Acerbi and Nicola Fusco.
\newblock A regularity theorem for minimizers of quasiconvex integrals.
\newblock {\em Arch. Rat. Mech. Anal.}, 99(261):261–281, 1987.

\bibitem{A12}
H.~W. Alt.
\newblock {\em Lineare {F}unktionalanalysis}.
\newblock Springer-Verlag GmbH, April 2012.

\bibitem{BallOP}
John~M. Ball.
\newblock {\em {Some Open Problems in Elasticity}}, chapter I.1, pages 3--59.
\newblock Springer-Verlag, NY, 2002.

\bibitem{BOP91MS}
P.~Bauman, N.~C. Owen, and D.~Phillips.
\newblock {Maximal smoothness of solutions to certain
  Euler{\textendash}Lagrange equations from nonlinear elasticity}.
\newblock {\em Proc. R. Soc. Ed. Sec. A.}, 119(3-4):241--263, 1991.

\bibitem{BOP91}
P.~Bauman, N.~C. Owen, and D.~Phillips.
\newblock {Maximum Principles and a priori estimates for a class of problems
  from nonlinear elasticity}.
\newblock {\em Annales de l{\textquotesingle}Institut Henri Poincare (C) Non
  Linear Analysis}, 8(2):119--157, mar 1991.

\bibitem{BOP92}
P.~Bauman, N.~C. Owen, and D.~Phillips.
\newblock {Maximum Principles and a priori estimates for an incompressible
  material in nonlinear elasticity}.
\newblock {\em Communications in Partial Differential Equations},
  17(7):1185--1212, 1992.

\bibitem{BK19}
J.~Bevan and S.~Käbisch.
\newblock {Twists and shear maps in nonlinear elasticity: explicit solutions
  and vanishing Jacobians}.
\newblock {\em Proc. R. Soc. Ed. Sec. A.}, 2019.

\bibitem{BY07}
J.~Bevan and X.~Yan.
\newblock Minimizers with topological singularities in two dimensional
  elasticity.
\newblock {\em {ESAIM}: Control, Optimisation and Calculus of Variations},
  14(1):192--209, sep 2007.

\bibitem{JB14}
Jonathan~J. Bevan.
\newblock On double-covering stationary points of a constrained {D}irichlet
  energy.
\newblock {\em Annales de l{\textquotesingle}Institut Henri Poincare (C) Non
  Linear Analysis}, 31(2):391--411, 2014.

\bibitem{JB17}
Jonathan~J. Bevan.
\newblock {A condition for the {H}ölder regularity of local minimizers of a
  nonlinear elastic energy in two dimensions}.
\newblock {\em Arch. Rat. Mech. Anal.}, 225(1):249--285, 2017.

\bibitem{BeDe21}
Jonathan~J. Bevan and Jonathan H.~B. Deane.
\newblock {Energy minimizing N-covering maps in two dimensions}.
\newblock {\em Calc. Var. Par. Dif. Eq.}, 60(4), 2021.

\bibitem{D1}
Jonathan~J. Bevan and Marcel Dengler.
\newblock A uniqueness criterion and a counterexample to regularity in an
  incompressible variational problem.
\newblock {\em arxiv:2205.07694}, 2022.

\bibitem{CY07}
Sungwon Cho and Xiaodong Yan.
\newblock {On the singular set for Lipschitzian critical points of polyconvex
  functionals}.
\newblock {\em Journal of Mathematical Analysis and Applications},
  336(1):372--398, 2007.

\bibitem{CG20}
Sergio Conti and Franz Gmeineder.
\newblock {Quasiconvexity and Partial Regularity, arXiv:2009.13820, 2020}.

\bibitem{C14}
Judith~Campos Cordero.
\newblock {\em {Regularity and Uniqueness in the Calculus of Variations}}.
\newblock PhD thesis, Oxford University, 2014.

\bibitem{CLM17}
Giovanni Cupini, Francesco Leonetti, and Elvira Mascolo.
\newblock {Local Boundedness for Minimizers of Some Polyconvex Integrals}.
\newblock {\em Arch. Ration. Mech. Anal.}, 224:269–289, 2017.

\bibitem{D2}
Marcel Dengler.
\newblock Everywhere regularity results for a polyconvex functional in finite
  elasticity.
\newblock {\em arXiv:2205.07694}, 2022.

\bibitem{EM01}
Luca Esposito and Giuseppe Mingione.
\newblock {Partial Regularity for Minimizers of Degenerate Polyconvex
  Energies}.
\newblock {\em Journal of Convex Analysis}, 8(1):1–38., 2001.

\bibitem{Ev86}
L.~C. Evans.
\newblock {Quasiconvexity and partial regularity in the Calculus of
  Variations}.
\newblock {\em Arch. Rat. Mech. Anal.}, 95(3), 1986.

\bibitem{LE10}
L.~C. Evans.
\newblock {\em {Partial Differential Equations: Second Edition (Graduate
  Studies in Mathematics)}}.
\newblock American Mathematical Society, 2010.

\bibitem{EVGA99}
L.C. Evans and R.F. Gariepy.
\newblock On the partial regularity of energy-minimizing, area-preserving maps.
\newblock {\em Calculus of Variations and Partial Differential Equations},
  9(4):357--372, dec 1999.

\bibitem{FR95}
Martin Fuchs and Jürgen Reuling.
\newblock {Partial regularity for certain classes for polyconvex functionals
  related to nonlinear elasticity}.
\newblock {\em Manuscripta mathematica}, 87(1):13--26, 1995.

\bibitem{FH94}
Nicola Fusco and John~E. Hutchinson.
\newblock {Partial regularity and everywhere continuity for a model problem
  from non-linear elasticity}.
\newblock {\em Austral. Math. Soc. (Series A)}, 57:158--169., 1994.

\bibitem{K11}
Aram~L. Karakhanyan.
\newblock Sufficient conditions for regularity of area-preserving deformations.
\newblock {\em Manuscripta Mathematica}, 138(3-4):463--476, 2011.

\bibitem{K13}
Aram~L. Karakhanyan.
\newblock Regularity for energy-minimizing area-preserving deformations.
\newblock {\em Journal of Elasticity}, 114(2):213--223, 2013.

\bibitem{KM07}
Jan Kristensen and Giuseppe Mingione.
\newblock {The Singular Set of Lipschitzian Minima of Multiple Integrals}.
\newblock {\em Arch. Rat. Mech. Anal.}, page 341–369, 2007.

\bibitem{KT03}
Jan Kristensen and Ali Taheri.
\newblock {Partial Regularity of Strong Local Minimizers in the
  Multi-Dimensional Calculus of Variations}.
\newblock {\em Arch. Rat. Mech. Anal.}, 170(1):63--89, nov 2003.

\bibitem{MS03}
S.~Müller and V.~{\v{S}}ver{\'{a}}k.
\newblock Convex integration for {L}ipschitz mappings and counterexamples to
  regularity.
\newblock {\em Annals of Mathematics}, 157(3):715--742, may 2003.

\bibitem{SS09}
Sabine Schemm and Thomas Schmidt.
\newblock {Partial regularity of strong local minimizers of quasiconvex
  integrals with $(p, q)-$growth}.
\newblock {\em Pro. Roy. Soc. Ed. Sec. A Math.}, 139(3):595--621, 2009.

\bibitem{Sz04}
L.~Szekelyhidi, Jr.
\newblock {The Regularity of Critical Points of Polyconvex Functionals}.
\newblock {\em Arch. Rat. Mech. Anal.}, 172(1):133--152, apr 2004.

\bibitem{GT12}
G.~Teschl.
\newblock {\em {Ordinary Differential Equations and Dynamical Systems (Graduate
  Studies in Mathematics)}}.
\newblock American Mathematical Society, 2012.

\bibitem{Y06}
X.~Yan.
\newblock {Maximal Smoothness for solutions to equilibrium equations in 2D
  nonlinear elasticity}.
\newblock {\em Proc. Amer. Math. Soc.}, 135(6):1717–1724, 2006.

\end{thebibliography}
\bibliographystyle{plain}
\vspace{0.5cm}
\textsc{Acknowledgements:} The author is appreciative to the Department of Mathematics at the University of Surrey and was funded by the Engineering \& Physical Sciences Research Council (EPRSC). Thanks to Jonathan J. Bevan and Bin Cheng for comments and discussions.

\end{document}